\documentclass[11pt]{article}

\usepackage{amssymb,amsmath,amsthm,amsfonts}

\usepackage{epsfig}

\usepackage{rotating}

\usepackage{subfigure}

\def\R{{\mathbb{R}}}
\def\T{{\mathbb{T}}}

\begin{document}

\title{Application of the $t$-model of optimal prediction to the estimation of
the rate of decay of solutions of the Euler equations in two and three dimensions}
\author{Ole H. Hald and Panagiotis Stinis\\ 
\\
Department of Mathematics \\
    University of California \\
and \\
Lawrence Berkeley National Laboratory \\
    Berkeley, CA 94720}
\date {}

\maketitle

\begin{abstract}
The ``$t$-model"  for dimensional reduction is applied to the estimation of the rate of decay
of solutions of the Burgers equation and of the Euler equations in two and three
space dimensions. The model was first derived in a statistical mechanics context,
but here we analyze it purely as a numerical tool and prove its convergence. 
In the Burgers case the model captures the rate of decay exactly, as was already previously 
shown. 
For the Euler equations in two space dimensions, the model preserves energy as it should. In three dimensions, we find a power law decay in time and observe a temporal intermittency.
\end{abstract}

\section{Introduction}
Despite the rapid increase in available computational power there are still many systems which cannot be studied numerically without prior simplification. In earlier work \cite{CHK00,CHK3}, we and others have derived methods for reducing
the number of variables one has to solve for in complex problems, based
on statistical projections. A special case, a long memory model called the "$t$-model", was thought to be particularly applicable to problems
in fluid dynamics \cite{CS05,alder}, where temporal correlations decay slowly. An earlier application \cite{bernstein} of the $t$-model to the
estimation of the rate of decay of solutions of the Burgers equation
yielded remarkably accurate results. 

In the present paper we use the $t$-model equations to reduce the number of variables in spectral methods 
and prove its convergence as the number of Fourier components increases . We then apply it to
the estimation of the rate of decay of solutions of Euler's equations in
two and three space dimensions.  We do not address the claim implicit in earlier
work, that the $t$-model equations may yield acceptable results even
when the number of variables remains finite. The results
we obtain are, however, surprisingly accurate, and a full analysis
may well have to go through some version of the arguments on the basis of
which the $t$-model was originally derived. 
Note that unlike previous damping methods for allowing spectral calculations
to proceed to significant time spans
(e.g. \cite{foias,moser,pasquetti,piomelli,scotti,she,smith}, the $t$-model
equations contain no adjustable parameters and
is guaranteed to remain stable.

The paper is organized as follows. In Section 2.1 we present the derivation of the $t$-model. In Section 2.2 we prove some results about its behavior for systems that conserve the $L_2$ norm of the solution and construct numerical methods that respect these properties. In Section 3, the $t$-model for the 3D Euler equations is constructed. In Section 4, we apply the $t$-model to the 1D inviscid Burgers equation and the 2D and 3D Euler equations and discuss how the numerical 
results compare to the known theoretical results.


\section{The $t$-model}

We begin with a system of ordinary differential equations
\begin{eqnarray}
\frac{d}{dt}v(t) &=& f\left(v(t),w(t)\right) , \, v(0)=x\\
\frac{d}{dt}w(t) &=& g\left(v(t),w(t)\right), \, w(0)=y 
\end{eqnarray}
Here $x,v \in \mathbb{R}^n$,
$y,w\in\R^m$
and $f:\R^n\times\R^m\rightarrow\R^n$,
$g:\R^n\times\R^m\rightarrow\R^m$,
and $t$ is time.
We think of  $v$ as the slow (resolved) variables and of $w$ as the 
fast (unresolved) variables.

We assume that the system (1)--(2) conserves energy and that the energy 
is given by
\begin{equation}
E=\frac{1}{2}(\|v\|^2+\|w\|^2).
\end{equation}

Here $\|v\|$ and
$\|w\|$ are the norms corresponding to the inner products
$(v,v')=\sum_{i=1}^n v_i v_i '$ and 
$(w,w')=\sum_{i=1}^m w_i w_i '.$
It follows from the conservation of energy that
\begin{eqnarray}
\left(v,f(v,0)\right) &=& 0 \\
\| f(v,0)\|^2+ \left(v,Df_v(v,0)\cdot f(v,0)\right)&=&0 \\
\| g(v,0)\|^2+ \left(v,Df_w(v,0)\cdot g(v,0)\right)&=&0
\end{eqnarray}
for all $v\in \R^n$.
Here the $n\times m$ matrix $D_w f(v,0)$ consists of the derivatives of
$f(v,w)$, evaluated at $w=0$.

To establish (4) we differentiate both sides of (3) with respect to $t$ and use 
(1),(2).
This gives
\begin{equation}
\frac{d}{dt} E=\left(v,f(v,w)\right)+\left(w,g(v,w)\right)=0.
\end{equation}
Since $v,w$ can be given any initial values we see that (7) is an 
identity in 
$v,w$.
In particular it holds when $w=0$ so $\left(v,f(v,0)\right)=0$ for all 
$v$.

To prove (5) we use a variational argument.
Let $a\in\R^n$.
Since (4) remains true when we replace $v$ by $v+\epsilon a$
it follows from Taylor's formula that
\begin{eqnarray*}
0&=&\left(v+\epsilon a,f(v+\epsilon a,0)\right) \\
&=&\left(v,f(v,0)\right)+\epsilon\left[\left(a,f(v,0)\right)+(v,D_v 
f(v,0)\cdot a) \right]+ O(\epsilon^2).
\end{eqnarray*}
But $\left(v,f(v,0)\right)=0$ so dividing by $\epsilon$,
letting $\epsilon\rightarrow0$ 
and setting $a=f(v,0)$ yields (5).

The proof of (6) is similar. 
Let $w=\epsilon \, b=\epsilon g(v,0)$.
Using Taylor's formula in (7) we obtain
\begin{eqnarray*}
0&=&\left(v,f(v,\epsilon b)\right)+\left(\epsilon b,g(v,\epsilon b)\right ) \\
&=&\left(v,f(v,0)\right)+\epsilon \left[\left(v,D_w f(v,0)\cdot b\right)+
\left(b,g(v,0)\right)\right]+O(\epsilon^2).
\end{eqnarray*}
To get (6) we divide by $\epsilon$ and let
$\epsilon\rightarrow 0$.

\subsection{Derivation of the $t$-model}
In this section we will derive and analyze approximations for systems 
which conserve energy and which can be written as (1)--(2).

Let $y=0$.
Since the energy is conserved, $w(t)=O(t)$.
Expanding $f(v,w)$ around $w=0$ we see that
\begin{equation}
\frac{d}{dt}v(t)=f\left(v(t),0\right)+D_w 
f\left(v(t),0\right)\cdot w(t)+O(t^2)
\end{equation}
with $y=0,$ equation (2) gives
\begin{eqnarray*}
w(t)&=&\int_0^t g\left(v(\tau),w(\tau)\right)d\tau\\
&=&\int_0^t g\left(v(\tau),0\right)d\tau+\int_0^t O(\tau)d\tau\\
&=&\int_0^t 
g\left(v(t),0\right)d\tau+\int_0^t O(t-\tau)d\tau+O(t^2)\\
&=&t g\left(v(t),0\right)+O(t^2)
\end{eqnarray*}
Inserting the expression  for $w$ in (8) and disregarding the 
$O(t^2)$
terms we arrive at the $t$-model
\begin{equation}
\frac{d}{dt}v(t)=f\left(v(t),0\right)+t D_w f\left(v(t),0\right)\cdot
g\left(v(t),0\right).
\end{equation}
It is called   the $t$-model because it contains the factor $t$, and 
because
it can be derived---for the cases we are interested in---as the zero
variance limit of the $t$-damping equations studied by Chorin, Hald 
and Kupferman \cite{CHK00}.

\subsection{Properties of the $t$-model and associated numerical methods}

The energy for a solution of the $t$-model is not constant, but 
decreases.
Indeed,  it follows from (4),(6),(9) that
\begin{eqnarray}
\frac{d}{dt}\frac{1}{2}\|v\|^2&=&\left(v,f(v,0)\right)+t\left(v,D_w 
f(v,0)\cdot
g(v,0)\right) \notag \\
&=&-t\|g(v,0)\|^2.
\end{eqnarray}
Thus the last term in (9) acts as a (non-linear) viscosity term.
Similar results have been obtained for the $t$-damping method applied 
to 
Hamiltonian systems, see \cite{CHK3}.

To solve eq.(9) we look for numerical methods where the energy 
decreases in each time step.
Let $F\left(v(t),t\right)$ denote the right-hand side of eq.(9) and
consider Runge-Kutta methods of the form
\begin{eqnarray*}
k_i&=&F\left(v^n+h\sum_{j=1}^s a_{ij} k_j, t_n+hc_i\right)\\
v^{n+1}&=&v^n+h\sum_{i=1}^s b_i k_i
\end{eqnarray*}
with $t_n=nh,$ $n=0,1,\ldots$.
Set $V_i=v^n+h\sum_{j=1}^s a_{ij} k_j$ for $i=1,\ldots,s.$

\vskip18pt

\noindent
{\bf Theorem}
Let $b_i a_{ij} + b_j a_{ji}-b_i b_j=0$ for $i,j=1,\ldots,s$ and assume 
that
\begin{equation}
\sum_{i=1}^s b_i=1\ \ ,\ \ \sum_{i=1}^s b_i c_i=\frac{1}{2}
\end{equation}
with $b_i,c_i\geq0$.
There is a $\tilde V$ in the convex hull of $V_1,\ldots,v_s$ such that
$$\frac{1}{2}\|v^{n+1}\|^2-\frac{1}{2}\|v^n\|^2=-h\left(t_n+\frac{h}{2}\right)\|g(\tilde 
V,0)\|^2.
$$

\noindent
{\bf Remark}
A numerical method that satisfies the assumptions in the theorem will 
be 
symplectic and at least second order.
The simplest example is the implicit midpoint rule.
It has $s=b_1=1$ and $a_{11}=c_1=\frac{1}{2}$.
Methods of higher order (4,5,6,8) can be found in \cite{hair} [p.207, p.209, p317].

\noindent
{\bf Proof}
We begin by expanding $\|v^{n+1}\|^2$.
After adding and subtracting $h\sum_j a_{ij}k_j$ we get
\begin{eqnarray*}
(v^{n+1},v^{n+1})= (v^n,v^n)&+&h\sum_i b_i(k_i,v^n+h\sum_j a_{ij}k_j)\\
&+&h\sum_j b_j(v^n+h\sum_i a_{ji}k_i,k_j)\\
&-&h^2\sum_{ij}[b_i a_{ij}+b_j a_{ji}-b_i b_j](k_i,k_j).
\end{eqnarray*}
The last sum vanishes.
Using the definitions of $v_i$ and $k_i$ yields
\begin{eqnarray*}
\|v^{n+1}\|^2-\|v^n\|^2&=&h\sum_i b_i\left(F(V_i,t_n+hc_i),V_i\right)\\
&+&h\sum_j b_j\left(V_j,F(V_j, t_n+hc_j)\right).
\end{eqnarray*}
Now $V_i$ are real, so the two sums are equal.
Consequently
\begin{eqnarray*}
\|v^{n+1}\|^2-\|v^n\|^2&=&2h\sum_j 
b_j\left(V_j,f(V_j,0)+(t_n+jc_j)D_wf(V_j,0)\cdot g(V_j,0)\right)\\
&=&-2h\sum_j b_j(t_n+hc_j)\|g(V_j,0)\|^2
\end{eqnarray*}
where we have used (4),(6).
Let
$$G_0=\min\|g(V_j,0)\|^2=\|g(V_{j_0},0)\|^2\ $$
$$G_1=\max\|g(V_j,0)\|^2=\|g(V_{j_1},0)\|^2.$$
Using (11) we conclude that
$$-2h\left(t_n+\frac{h}{2}\right)G_1\leq\|v^{n+1}\|^2-\|v^n\|^2\leq 
-2h\left(t_n+\frac{h}{2}\right)G_0.$$
Set $\tilde V=\theta V_j+(1-\theta)V_j$.
Since $g$ $(V,0)$ is continuous there is a $\Theta\in[0,1]$ such that
$$\|v^{n+1}\|^2-\|v^n\|^2=-2h\left(t_n+\frac{h}{2}\right)\|g(\tilde 
V,0)\|^2.$$
This completes the proof.


\section{The $t$-model for the Euler equations}

The Euler equations describe the flow of an incompressible, inviscid 
fluid
in two or three dimensions.
Here we look at flows in a cube with periodic boundary conditions and 
consider
two kinds of approximations.
First we use the Fourier method to obtain approximate solutions of 
Euler's
equations.
This leads to a large system of ordinary differential equations.
Secondly we use the $t$-model to reduce the number of variables.
The questions are: Does this method converge and what
are the numerical results?

The Euler equations for 3 dimensional flows are
\begin{eqnarray}
\partial_t u+(u\cdot\nabla)u&=&-\nabla p\\
\nabla \cdot u&=&0
\end{eqnarray}
where $u=u(x,t)$ is the velocity, $p$ is the pressure and $t$  is time.
We look for $2\pi$ periodic solutions and use the notation
$u\cdot\nabla=\sum_j u^j\partial_j$ where $u=(u^1,u^2,u^3)^T$, $T$ 
means
transpose and $\partial_j u=\partial u/\partial x_j$ for $j=1,2,3$.
By taking the divergence of (12) and using (13) we see that
$$-\Delta p=\nabla\cdot(u\cdot\nabla)u$$
where $\Delta=\sum_j\partial_j^2$.
Set $\int p=0$.
We can then solve for $p$ and conclude from (12) that
\begin{equation}
\partial_t 
u=-\left[(u\cdot\nabla)u+\nabla\cdot(-\Delta)^{-1}\nabla\cdot(u\cdot\nabla)u\right].
\end{equation}
Next we expand $u$ in a Fourier series
\begin{equation}
u(x,t)=\sum_k u_k(t) e^{ikx}.
\end{equation}
It follows from (13) that $k\cdot u_k\equiv0$ for $k\neq0$.
By combining (12) and (15) and proceeding formally we get for 
$k\neq0$

\begin{equation}
\frac{d}{dt}u_k=-i\sum_{p+q=k} k\cdot u_p A_k u_q
\end{equation}
where $A_k=I-kk^T/|k|^2.$ (For more details, see e.g. \cite{canuto}).
Since $u$ is real, $u_{-k}=\bar u_k$.
Finally we assume that $u_0\equiv0$.
As $\int u$ is constant in time this amounts to a restriction of the 
initial 
data.

We get the Fourier method by setting $u_k\equiv0$ for $|k|_{\infty}>m $ and for all time 
and setting $k\cdot u_k=0$ and $u_{-k}=\bar u_k$ for $t=0$.
In this way (16) becomes a closed system of ordinary differential equations that conserves energy 
and
yields solutions that are incompresible.
To express (16) in the form (1)--(2) we set $v_k=u_k$ if $|k|_{\infty}\leq n,$
and $w_k=u_k$ if $n<|k|_{\infty} \leq m$.
Let $F=\{|k|_{\infty}\leq n\}$ and $G=\{n<|k|_{\infty}\leq m \}$.
If $k\in F$ then
\begin{eqnarray*}
 \frac{d v_k}{dt}=- i \biggl[ \underset{p \in F ,\, q \in F }{\underset{p+q=k  }{ \sum}} k \cdot v_{p} 
A_{k} v_{q} + \underset{p \in F ,\, q \in G}{\underset{p+q=k  }{ \sum}} k \cdot v_{p} 
A_{k} w_{q} \\
\underset{p \in G ,\, q \in F }{\underset{p+q=k  }{ \sum}} k \cdot w_{p} 
A_{k} v_{q} +\underset{p \in G  ,\, q \in G}{\underset{p+q=k  }{ \sum}} k \cdot w_{p} 
A_{k} w_{q} \biggr]. 
\end{eqnarray*}
The equation for $dw_k/dt$ is the same, except that $k\in G$.
Following the derivation of (9) we obtain the $t$-model for 
Euler's equations
\begin{eqnarray}
\frac{d}{dt}v_k=- i  \biggl[ \underset{p \in F ,\, q \in F }{\underset{p+q=k  }{ \sum}} k \cdot v_{p} 
A_{k} v_{q} +  \underset{p \in F  ,\, q \in G}{\underset{p+q=k  }{ \sum}} k \cdot v_{p} 
A_{k} (-i) t  \underset{r \in F ,\,  s \in F }{\underset{r+s=q }{ \sum}} q \cdot v_{r} 
A_{q} v_{s} \\
+ \underset{p \in G ,\, q \in F }{\underset{p+q=k  }{ \sum}} k \cdot 
(-i) t  \underset{r \in F  ,\, s \in F }{\underset{r+s=p }{ \sum}} p \cdot v_{r} 
A_{p} v_{s} A_{k} v_{q} \biggr] \notag
\end{eqnarray}
If $k\cdot v_k=0$ and $v_{-k}=\bar v_k$ at $t=0$,
it holds for all time.
Moreover, a direct calculation shows that the energy decays, i.e.
\begin{equation}
\frac{d}{dt}\frac{1}{2}\sum_{k\in F}|v_k|^2=-t\sum_{k\in G} | \underset{p \in F  ,\, q \in F }{\underset{p+q=k  }{ \sum}} k\cdot v_p A_k v_q | ^2.
\end{equation}
This is the complex analogue of (10).
Finally we compare the solution $u(x,t)$ of Euler's equations with the
solutions generated by the $t$-model.
Let $\T^d$ be the torus of length $2\pi$ in ${\R}^d$ with
$d=2,3$ and set
$$v(x,t)=\sum_{k\in F} v_k(t)e^{ik\cdot x}.$$

\noindent
{\bf Theorem}
Assume that $u(t)\in H^s(\T^d)$ for $s\geq3$, $d=2,3$ and $0\leq t\leq T$.
Then
$$\|u(t)-v(t)\|_{L^2(\T^d)}\leq\frac{{\rm const}(T)}{n^{s-1}}\cdot
\max_{0\leq t\leq T}\|u\|_{H^s(\T^d)}.$$

\noindent
{\bf Proof}
Will be presented elsewhere.

Thus if the solution of Euler's equations is smooth for $t\leq T$
then the energy is constant, and the energy of the $t$-model should
converge to the energy of Euler's equations.
In our numerical experiments this holds for $d=2$ and fails for $d=3$.
This suggests  that the solutions of the Euler equation lose smoothness
when $d=3$, and may develop singularities.


\section{Numerical results}

In this section we present numerical results of the application of the $t$-model to the 1D inviscid Burgers equation and the 2D and 3D Euler equations. The equations of motion for the Fourier modes were solved by a
Runge-Kutta-Fehlberg method  (\cite{hair}) with the tolerance set to $10^{-10}.$ 
Note that due to the quadratic nonlinearity and the form of the $t$-model term (see Eq. 17), the right hand side of the equation for each Fourier mode in the $t$-model contains interactions with Fourier modes of at most double the wavevector. So, the ratio of the number of unresolved modes in $G$ to the resolved modes in $F$ is 1.

The different terms appearing in the right hand side of the equations for the reduced model can be computed in 
real space using the Fast Fourier Transform (FFT). Since for a reduced model of size $N$ in each spatial direction we include $N$ unresolved modes in each direction, the arrays involved in the FFTs should be of size $2N$. The fact that the $t$-model term can be computed using the FFT makes the numerical 
implementation of the $t$-model computationally efficient. Moreover, the FFT calculations involved 
are dealiased by construction and thus no extra (e.g. $3/2$ rule) dealiasing is needed. For a calculation involving N modes in each direction, i.e. $N/2$ positive and $N/2$ negative, we perform FFTs of size $2N$, i.e. $N$ positive and $N$ negative modes. But we are 
interested only on the right hand side of the equations for the first $N/2$ modes. This means (see \cite{canuto}) that to avoid 
aliasing (in the $t$-model term) we need for the total number of modes M used to satisfy the following inequality: $-N-N/2 \geq N/2-1-M$ which yields $M \geq 2N.$ But $2N$ is exactly how many modes 
we use in the FFTs, and thus the $t$-model term calculation through FFTs is dealiased by construction.

In order to study the 
asymptotic decay rate of the energy in the resolved modes, one has to evolve the system for long times. We evolved each case up to time t=100, so that we have enough points to perform an accurate estimate of the decay rate exponent. The need to perform calculations for long times, prevented us from conducting 
numerical experiments of larger size (within reasonable time) for the 3D Euler equations on a single processor workstation. However, the fact that the $t$-model term can be computed using the FFT means that the implementation 
of the model in any existing parallel spectral Navier-Stokes code is straightforward and we expect to report results of such simulations in the near future.

\begin{figure}
\centering
\subfigure[]{\epsfig{file=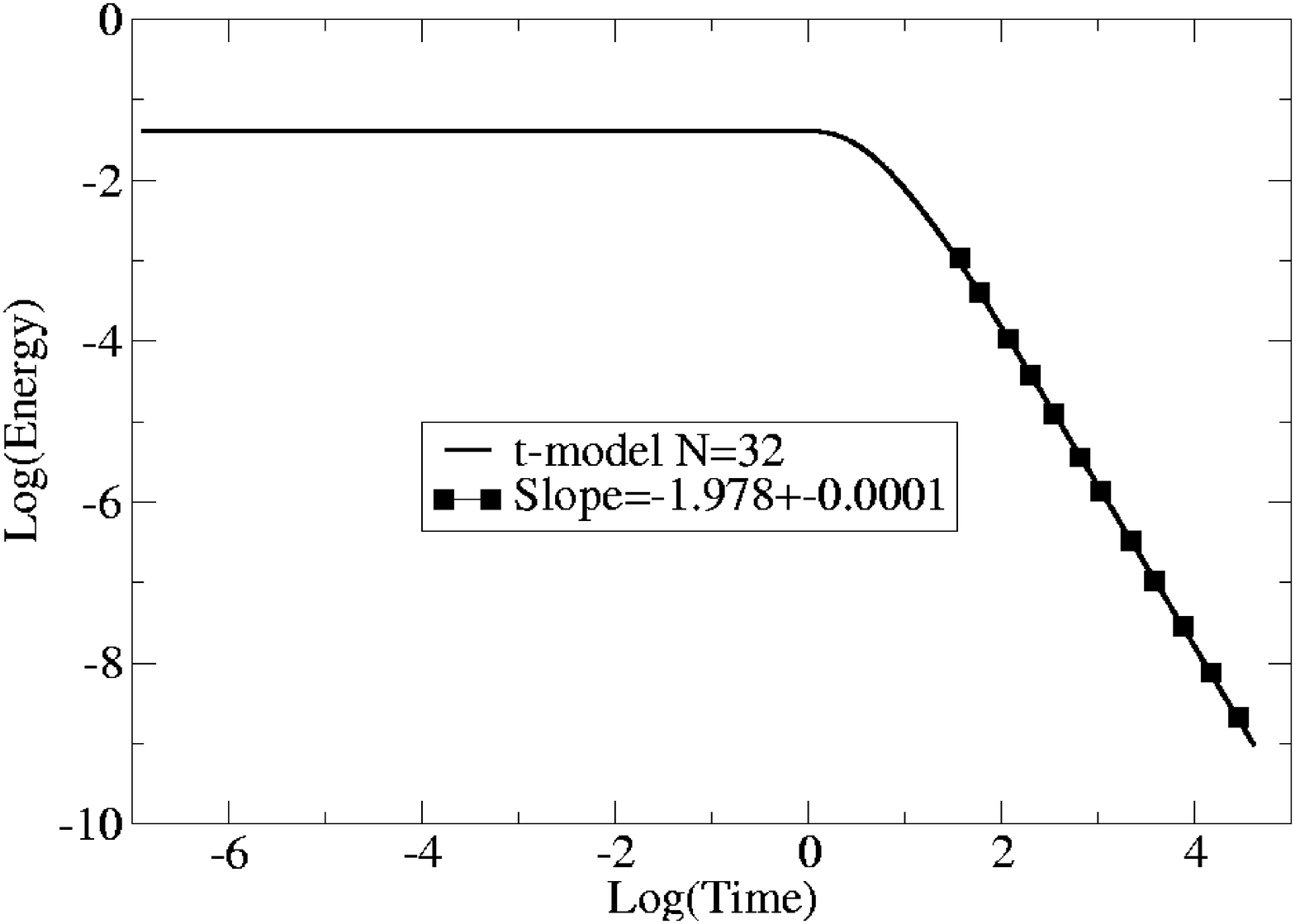, width=2.in}}
\qquad
\subfigure[]{\epsfig{file=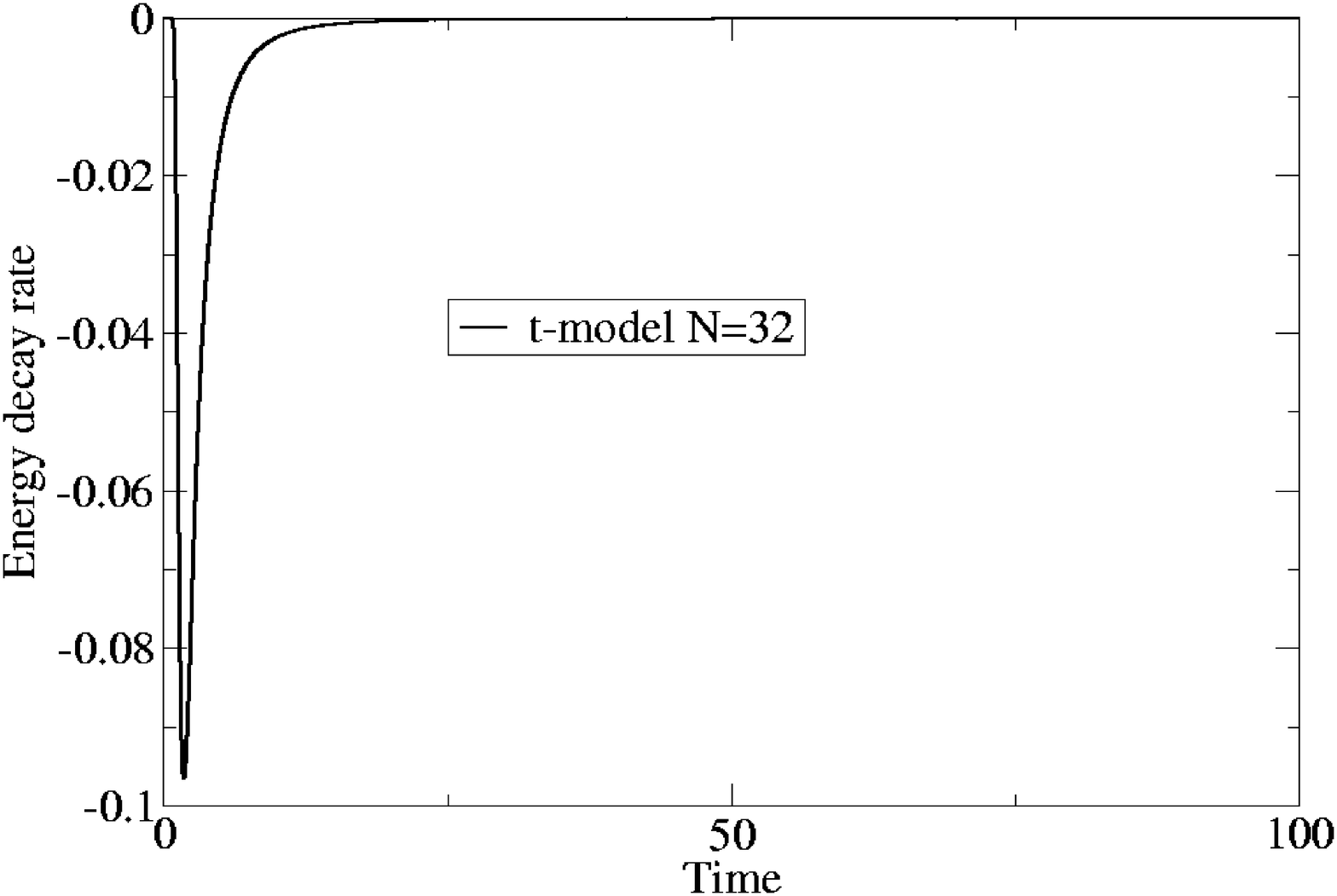,width=2.in}}
\caption{(a) Energy evolution of the $t$-model with $N=32$ modes for the inviscid Burgers equation. (b) Evolution of the energy decay rate. }
\label{fig_burgers}
\end{figure}

\begin{figure}
\centering
\subfigure[]{\epsfig{file=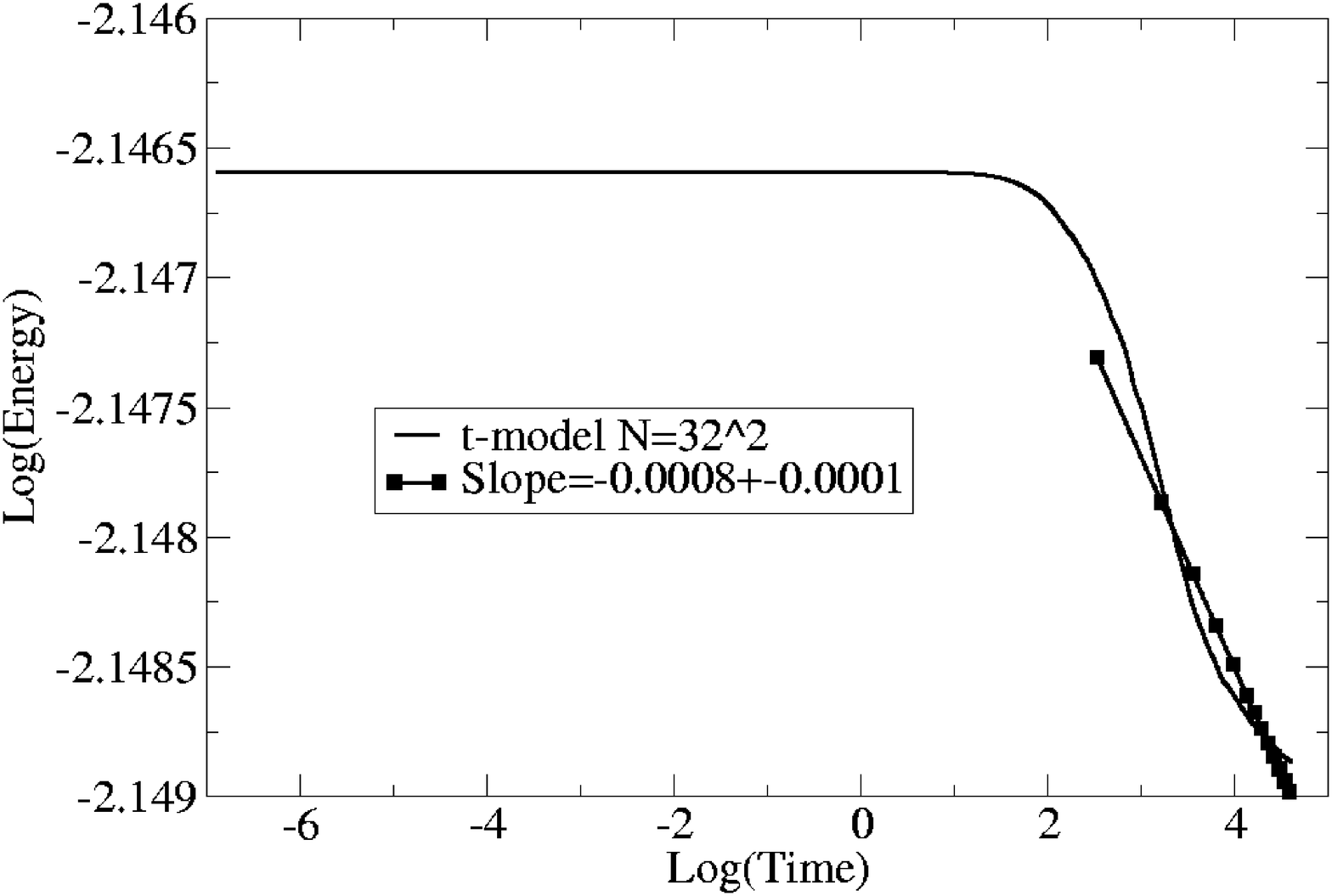, width=2.in}}
\qquad
\subfigure[]{\epsfig{file=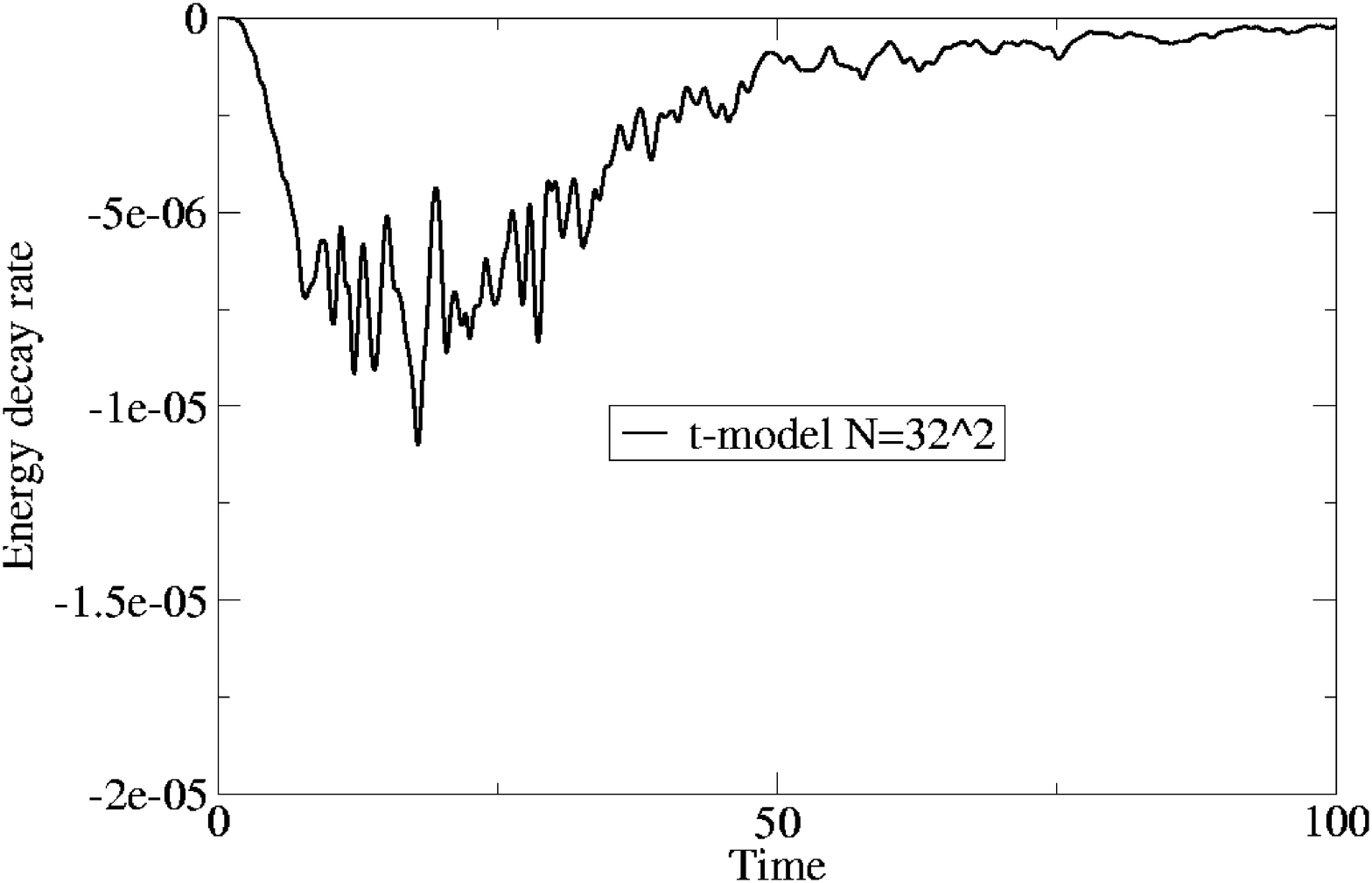,width=2.in}}
\caption{(a) Energy evolution of the $t$-model with $N=32^2$ modes for the 2D Euler equations (note the extremely slow decay reflected in the very small slope of the linear fit). (b) Evolution of the energy decay rate. Again, note the very small values of the energy decay rate. }
\label{fig_euler2d}
\end{figure}

\begin{figure}
\centering
\subfigure[]{\epsfig{file=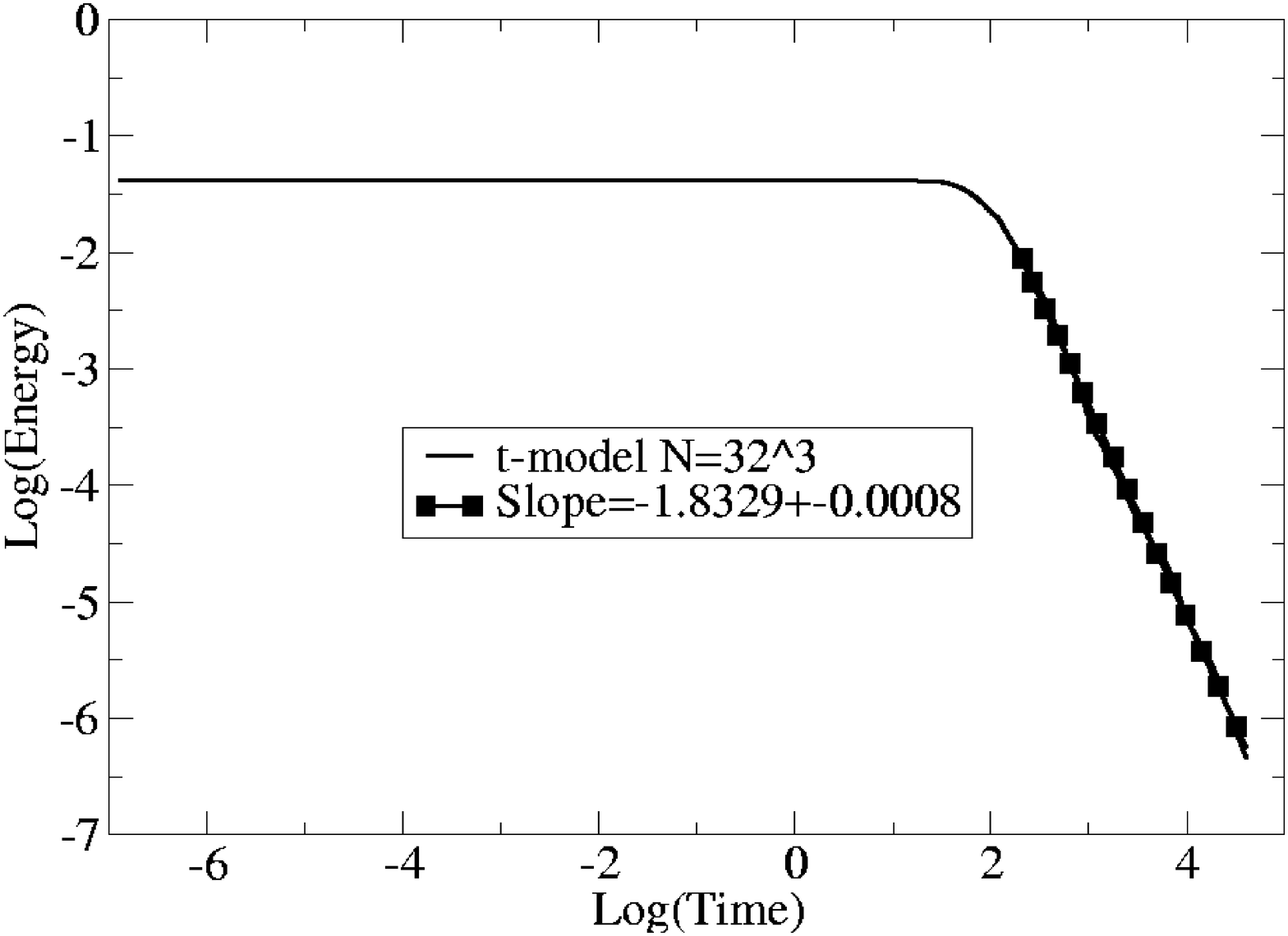, width=2.in}}
\qquad
\subfigure[]{\epsfig{file=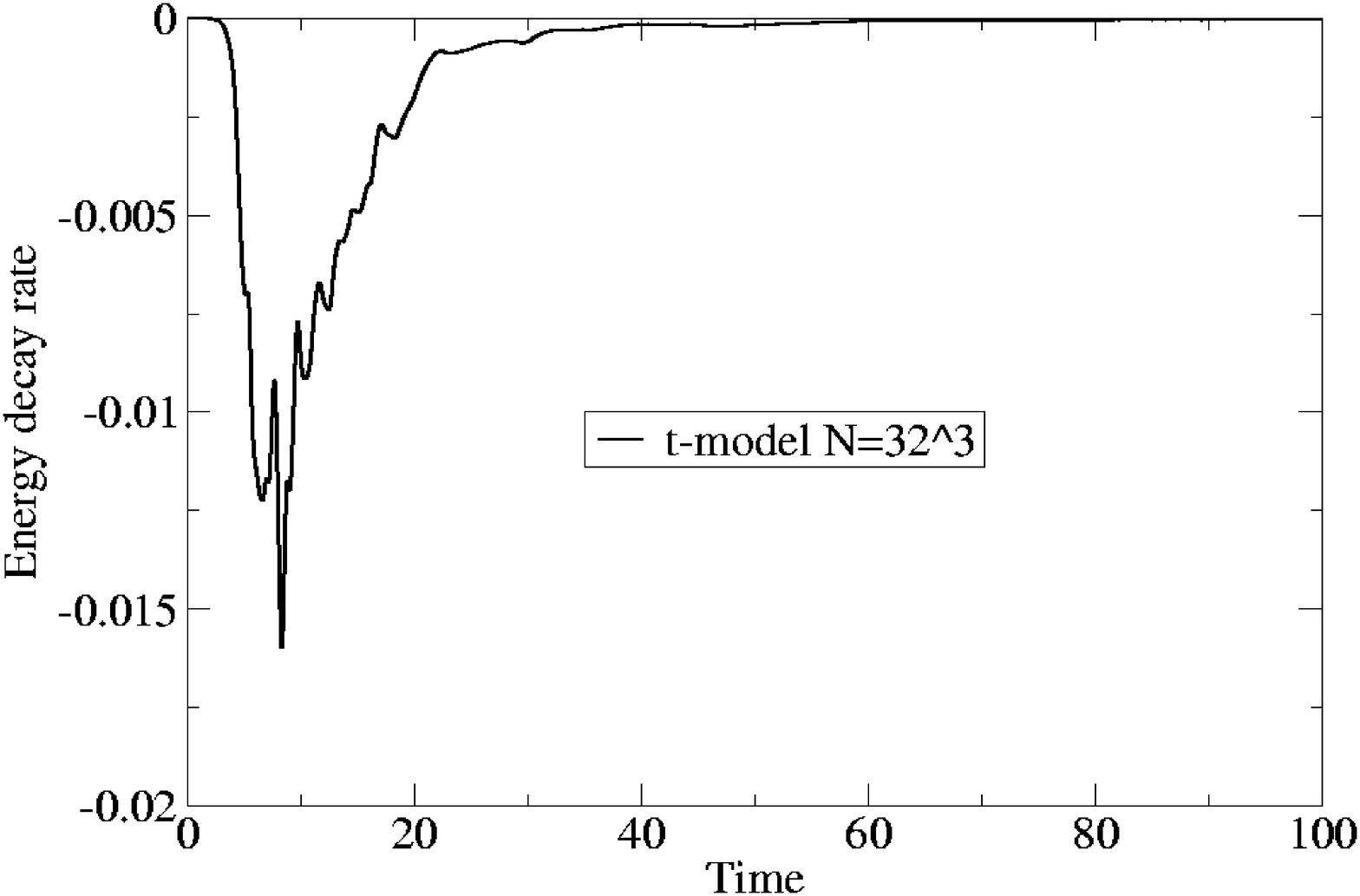,width=2.in}}
\caption{(a) Energy evolution of the $t$-model with $N=32^3$ modes for the 3D Euler equations. (b) Evolution of the energy decay rate. }
\label{fig_euler3d}
\end{figure}

Figures \ref{fig_burgers}, \ref{fig_euler2d} and \ref{fig_euler3d} present results of the application of the 
$t$-model to the 1D inviscid Burgers equation, the 2D Euler equations and the 3D Euler equations respectively. We present results for the evolution with time of the energy in the resolved modes and for the rate of energy decay (see Eq. 18). The numerical experiments are for resolved sets $F$ (and the corresponding sets $G$) of size $N=32,$ $N=32^2,$ and $N=32^3$ for the 1D, 2D and 3D cases respectively.

For the 1D inviscid Burgers equation $u_t + u u_x=0, $ the initial condition is $u_0(x)=\sin{x},$ which gives rise to a single shock wave at time $T=1.$ Until the moment of the formation of the shock wave the energy $E$ is practically constant and soon after the well known $t^{-2}$ energy decay law \cite{lax} is established. The energy decay rate shown in Figure \ref{fig_burgers}(b) reaches its peak around time $\tau=1.219$ after which it starts decreasing. The $t^{-2}$ energy decay regime is established soon after. The estimated exponent of the energy decay is estimated as $-1.9781\pm0.0001,$ using about 15000 points. It is interesting that the right energy decay law is captured with only $N=32$ Fourier modes.

For the 2D Euler equations the situation is drastically different. We present results for an incompressible, isotropic random initial condition with energy spectrum 
$E(k)=\exp(-2k)$ for the resolved modes and zero for the unresolved modes. After 100 units of time, the energy has decayed by $0.2\%.$ In other words, the smooth initial condition does not lose its smoothness. In fact, as one can see in 
Figure \ref{fig_euler2d}(a), after the insignificant energy decay, the energy establishes a plateau which signifies the absence of drain of energy out of the resolved range of modes. A linear fit (in log-log coordinates) of the energy evolution gives a slope of $-0.0008\pm 0.0001,$ where we used about 15000 points.

The qualitative difference between the 2D and 3D case is striking. For the 3D Euler equations we use as initial condition the Taylor-Green vortex given by
\begin{eqnarray*}
u^1(x,0)&=&\sin(x_1)\cos(x_2)\cos(x_3), \\
u^2(x,0)&=&-\cos(x_1)\sin(x_2)\cos(x_3),\\
u^3(x,0)&=&0
\end{eqnarray*} 
Note that the Taylor-Green initial condition is smooth having nonzero values only for the Fourier modes with $k_i=\pm1 ,\, i=1,2,3.$ A lot of numerical work (see e.g. \cite{gottlieb} and references therein) has been devoted to the investigation of whether the solutions of the 3D Euler equations with the Taylor-Green initial condition blow up in finite time. All calculations show a rapid increase in the value of the maximum vorticity, but are hampered by the fact that they run out of resolution around time $T=5.$ Note that the first peak of the energy decay rate that we find is around time $\tau=5.17.$

More interestingly, the decay of the energy appears to be organized in a collection of spikes of diminishing strength. This organization of the energy decay is reminiscent of the phenomenon of intermittency, i.e. bursts of activity followed by 
intervals of relative inaction on the part of the flow. Of course, the phenomenon of intermittency is not 
only of temporal nature, but has a spatial manifestation too. This is exhibited as concentration of the 
highest vorticity in small regions of the flow. The trend we observe in the decay of the energy seems to assign a specific purpose to the vorticity. Starting from a smooth initial condition, we have a steepening of the gradients in the field. This means that smaller scales are excited until the vorticity producing 
mechanism runs out of steam. Then we enter a period of relative inaction, until there is a restart of the 
mechanism of steepening. Energy is transferred again to the smaller scales and so forth. This scenario continues until there is no energy left in the large scales. After that, the flow just disintegrates and eventually comes to a halt. The purpose of vorticity mentioned above is to regulate the transfer of energy to the small scales \cite{C94}. This is reminiscent of the picture suggested by Moffatt, Kida and Okhitani \cite{moffatt} of the vortex structures  acting as the "sinews of turbulence" .

The slope of the linear fit (in log-log coordinates) for the energy decay is $-1.8329\pm0.0008$, where we have used about 15000 points. Estimates on the energy decay law for the 3D Euler equations are rare (see \cite{speziale} and references therein). Moreover, all the estimates concern the infinite space case while we use periodic boundary conditions. For the infinite space case, under the assumption of complete self-preservation, i.e. self-similarity for all scales from 0 to $\infty,$ one finds that the energy shoud decay as $t^{-1}$ \cite{speziale}. If the assumption of complete self-preservation is not satisfied, the energy is expected to decay as $t^{-\alpha}$ where $\alpha > 1.$ Note that the assumption of complete self-preservation is violated for the case of periodic boundary conditions and so the exponent of the energy decay should be larger than 1. However, it is not clear how the assumptions can be modified for this case. In 1D, the change of boundary conditions from infinite to periodic changed the exponent of energy decay from 1 to 2 \cite{lax}. It is not clear that this is also the case for 3D. If it is, then the numerical estimate -1.8329 for $\alpha$ becomes more plausible.


\section{Conclusions}
The problem of constructing reduced models for the Euler equations has been, and still is, a great 
challenge for scientific computing. The $t$-model proposed here should be considered a first step 
in deriving models directly from the equations without {\it ad hoc} approximations. It is based on numerical and physical observations about the behavior of the solution (more sophisticated reduced models for the Euler equations were constructed and simulated in \cite{S06}). Following \cite{bernstein}, we tested the model on the 1D inviscid Burgers equation for an initial condition that gives rise to a shockwave. The model captures the right time of formation of the shock and the right rate of decay of the energy of the solution. For the 2D Euler equations, the model preserves the energy as it should since the solution remains smooth for all times. The numerical results for the Taylor-Green vortex for the 3D Euler equations produce rates of decay compatible with current thinking, and suggest that the solution loses its smoothness in finite time.

The terms appearing in the reduced model can be efficiently implemented by the use of the FFT on appropriate arrays. This makes the incorporation of the model in existing pseudospectral algorithms rather straightforward. We plan to apply the model in a parallel setting which will allow a better assessment of the properties of the flow field that is predicted by the model.

\section{Acknowledgements}
We are grateful to Prof. G.I. Barenblatt, Prof. A.J. Chorin and Mr. J. Weare for many helpful discussions and comments. We are especially indebted to Dr. Yelena Shvets for her critical reading of the analysis and for moral support during times of bad results. This work was supported in part by the National Science 
Foundation under Grant DMS 04-32710, and by the Director,
Office of Science, Computational and Technology Research, U.S.\ Department of 
Energy under Contract No.\ DE-AC03-76SF000098.


\begin{thebibliography}{99}



\bibitem{CHK00}
Chorin, A.J., Hald, O.H. and Kupferman, R., Proc. Nat. Acad. Sc. USA 97 (2000) pp. 2968-2973.

\bibitem{CHK3}
Chorin, A.J., Hald, O.H. and Kupferman, R., Physica D 166 (2002) pp. 239-257.

\bibitem{CS05}
Chorin, A.J. and Stinis, P., Problem reduction, renormalization and memory,
Comm. App. Math. Comp. Sci. 1 (2005) pp. 1-27.


\bibitem{alder} Alder, B.  and Wainwright, T., Phys. Rev. A 1 (1970) pp. 1-12.


\bibitem{bernstein}
Bernstein, D., Multi. Mod. Sim. (2006) in press.



\bibitem{foias}
Foias, C. Holm, D.D., Titi, E.S., Physica D 152-153 (2001), 505-519.

\bibitem{moser}
Langford, J. and Moser, R., J. Fluid. Mech.
398 (1999) pp. 321-346.

\bibitem{pasquetti}
Pasquetti, R., J. Turb. 6 (2005) pp. 1-14.

\bibitem{piomelli}
Piomelli, U., Prog. Aero. Sci. 35 (1999) pp. 335-362.  

\bibitem{scotti}
Scotti, A. and Meneveau, C., Phys. Rev. Lett. 78 (1997) pp. 867-870.

\bibitem{she}
She, Z.S. and Jackson, E., Phys. Rev. Lett. 70 (1993) pp. 1225-1228.

\bibitem{smith}
Smith, L.M. and Woodruff, S.L., Ann. Rev. Fluid Mech. 30 (1998) pp. 275-310.


\bibitem{hair}
Hairer, E.,  N\"orsett, S.E.,  and Wanner, G., Solving Ordinary Differential Equations I, Springer, NY, 1987.

\bibitem{canuto}
Canuto, C., Hussaini, M.Y.,  Quarteroni, A.  and  Zang, T.A., Spectral Methods in Fluid Dynamics, Springer, NY, 1988.


\bibitem{lax}
Lax, P.D., Hyperbolic Systems of Conservation Laws and the Mathematical Theory of 
Shock Waves, SIAM Publications, Philadelphia, 1972.



\bibitem{gottlieb}
Don, W.S., Gottlieb, D., Shu, C.W. , Schilling, O. and Jameson, L., J. Sci. Comp. 24 (2005), pp. 569-595. 

\bibitem{C94}
Chorin, A.J.,  Vorticity and Turbulence, Springer, NY, 1994.


\bibitem{moffatt}
Moffatt, H.K, Kida, S. and Okhitani, K., J. Fluid Mech. 259 (1994) pp. 241-264.

\bibitem{speziale}
Speziale, C.G. and Bernard, P.S., J. Fluid Mech. 241 (1992) pp. 645-667.



\bibitem{S06}
Stinis, P., Technical Report LBNL-60899 (2006)  \& math.NA/0607108 (submitted to SIAM Multi. Mod. Sim.).


\end{thebibliography}
\end{document}